\documentclass[a4paper,english]{amsart}
\usepackage[T1]{fontenc}
\usepackage[latin1]{inputenc}
\usepackage{amsthm}
\usepackage{amssymb}
\usepackage{mathptmx}
\usepackage{babel}
\usepackage{bbm}
\usepackage{setspace}
\usepackage{microtype}

\numberwithin{equation}{section} %% Comment out for sequentially-numbered
\numberwithin{figure}{section} %% Comment out for sequentially-numbered
\theoremstyle{plain}
\newtheorem{thm}{Theorem}[section]
  \theoremstyle{plain}
  
  \theoremstyle{plain}
  
  \theoremstyle{plain}
    \newtheorem{prop}[thm]{Proposition}
      \theoremstyle{plain}
    \newtheorem*{fact*}{Fact}

\newcommand{\1}{\mathbbm{1}}
\newcommand{\maxx}{{\mbox{\it\footnotesize max}}}
\newcommand{\maxxx}{{\mbox{\it\tiny max}}}

\newcommand{\e}{\mathrm e}

\DeclareMathOperator{\diam}{diam}

\DeclareMathOperator{\card}{card}

\renewcommand{\hat}{\widehat}
\renewcommand{\phi}{\varphi}

\renewcommand{\tilde}{\widetilde}

\def\R{\mathbb{R}}
\def\H{\mathbb{D}}

\def\N{\mathbb{N}}
\def\Z{\mathbb{Z}}
\def\P{{\mathcal P}}

\begin{document}
 
\title[A note on the algebraic growth rate]{A note on the
algebraic growth rate of Poincaré series for Kleinian groups}

\author{Marc Kesseböhmer}
\address{Fachbereich 3 -- Mathematik und Informatik, Universität Bremen, Bibliothekstr.
1, D--28359 Bremen, Germany}
\email{mhk@math.uni-bremen.de}

\author{Bernd O. Stratmann}
\address{Mathematical Institute, University of St. Andrews, North Haugh, St.
Andrews KY16 9SS, Scotland}
\email{bos@st-andrews.ac.uk}

\subjclass[2000]{Primary 30F40; Secondary 37A05}
\date{October 26, 2009}
\keywords{Poincar\'e series, infinte ergodic theory, Kleinian groups}
 
 \begin{abstract}
In this note we employ infinite ergodic theory to 
derive  estimates for the algebraic growth rate of 
the Poincar\'e series
for a Kleinian group 
at its critical  exponent of convergence. 
\end{abstract}

 \dedicatory{Dedicated to S.J. Patterson on the occasion of his
 60th birthday.} 

\maketitle
\onehalfspacing
\section{Introduction and statements of result}

\noindent In this note we study the Poincar\'e series 
\[ \P(z,w,s) :=\sum_{g\in G}\e^{-s d(z,g(w))}  \]
of a geometrically finite, 
essentially free Kleinian group $G$ acting on 
the $(N+1)$-dimen\-sio\-nal hyperbolic space $\H,$ for arbitrary
$z,w\in\H$.  Here, $d(z,w)$ denotes the hyperbolic distance between 
$z$ and $w$, and $s \in \R$.
It is well known that a group of this type  is of $\delta$-divergence 
type, which means that $\P(z,w,s)$ diverges for $s$ equal to the exponent 
of convergence $\delta=\delta(G)$ of $\P(z,w,s)$.
We are in particular interested in the situation in which $G$ 
is a zonal group, that is, we always assume that $G$ has parabolic elements.  
For Kleinian groups of this type we then consider the partial 
Poincar\'e sum
\[ \P_{n}(z,w,s) :=\sum_{g\in G\atop |g|\leq n}\e^{-s d(z,g(w))},\]
for $n \in \N$, and where $| \, \cdot \, |$ denotes the word 
metric in $G$. The main result of this note is  the following
 asymptotic estimate for
these partial Poincar\'e sums, for $s$ equal to the exponent of 
convergence
$\delta$. Here, $r_{\maxx}$ denotes the 
maximal rank of the parabolic fixed points of $G$, and $\asymp$
denotes comparability, that is, $b_{n} \asymp c_{n}$ if and only if 
$(b_{n}/c_{n})$
is uniformly bounded away from zero and infinity, for two sequences 
$(b_{n})$ and $(c_{n})$ of positive 
real numbers.

\begin{thm}\label{main}
For a geometrically finite, essentially free, zonal Kleinian group 
$G$
and for each $z,w \in \H$,  we have  
\begin{eqnarray*}
    \P_{n}(z,w,\delta) & 
    \asymp & \left\{ \begin{array}{lcc}
    n^{ 2\delta-r_{\maxxx}} &
    \mbox{for } &   \delta<(r_{\maxx}+1)/2\\
     n/\log \, n & \mbox{for } & \delta= (r_{\maxx}+1)/2\\
    n & \mbox{for } & \delta>(r_{\maxx}+1)/2.\end{array}\right.\end{eqnarray*}
    \end{thm}  
 Note that the results in this note grew out of the authors closely 
  related recent 
  studies in \cite{KesseboehmerStratmann:07} of the so-called sum--level sets
  for regular continued fractions.  These sets are given by
\[  \mathcal{C}_{n}:=\{[a_{1},a_{2},\ldots] \in [0,1]:
\sum_{i=1}^{k} a_{i}=n, 
  \mbox{ for some } k\in\N\},\]
  where $[a_{1},a_{2},\ldots]$ denotes  the regular continued fraction expansion. 
  Inspired by a conjecture in \cite{FK}, it was shown  in
  \cite{KesseboehmerStratmann:07}
  that for the Lebesgue 
measure $\lambda(\mathcal{C}_{n})$ of these sets one
has that, with  $b_{n}\sim c_{n}$ denoting 
$\lim_{n\to\infty}b_{n}/c_{n}=1$,
\[
\lambda(\mathcal{C}_{n})\sim\frac{1}{\log_{2}n} \mbox{  and  } 
\sum_{k=1}^{n}\lambda\left(\mathcal{C}_{k}\right)
\sim\frac{n}{\log_{2} n}. 
\]
It is not hard to see that in here  the second asymptotic estimate  
implies Theorem \ref{main} for $G$ equal 
to the (subgroup of index $3$ of the) modular group $PSL_{2}(\Z)$.

\vspace{2mm}

{\it Acknowledgement.}
We would like to thank the {\em Mathematische Institut der Universit\"at G\"ottingen}  for the warm hospitality during our research visit 
in Summer 2009. In particular,  we would like to thank 
P. Mihailescu and S.J. Patterson for the excellent organisation of   the 
{\em International Conference: Patterson 60 ++},
which took place during the period of our visit.

   \section{Preliminaries}\subsection{The canonical Markov map}
   As already mentioned in the introduction, throughout, we
   exclusively consider a geometrically finite,
    essentially free, zonal Kleinian group $G$. 
    By definition (see \cite{KesseboehmerStratmann:04A}), a group of this type can be 
     written as a free product  $ G=H \ast \Gamma$,
   where $H=\langle h_{1},h_{1}^{-1}\rangle \ast \ldots 
   \ast \langle h_{u},h_{u}^{-1}\rangle$ denotes the 
     free product of finitely many elementary,  loxodromic 
     groups, and $\Gamma= \Gamma_{1} \ast 
     \ldots \ast  
     \Gamma_{v} $ denotes the 
     free product of finitely many  parabolic
     subgroups of $G$ such that $\Gamma_{i}
     =\langle 
     \gamma_{i_{1}},\gamma_{i_{1}}^{-1},\ldots, 
     \gamma_{i_{r_{i}}},\gamma_{i_{r_{i}}}^{-1} \rangle $ is the 
     parabolic subgroup of $G$ associated with the parabolic fixed 
     point $p_{i}$ of rank $r_{i}$. Clearly,  $\Gamma_{i} \cong {\Bbb 
     Z}^{r_{i}}$ and 
     $\gamma_{i_{j}}^{\pm}(p_{i})=p_{i}$, for all $j=1, \ldots , r_{i}$ 
     and $i=1,\ldots ,v$. Also, note that  $G$ 
     has no relations other than those which 
    originate from cusps of rank at least
    $2$, that is, those $\Gamma_{i}$ with $r_{i}>1$.  
    Without loss of generality,  we can assume that
    $G$ admits the choice of a 
    Poincar\'e polyhedron  $F$ with a 
    finite set $\mathcal{F}$ of faces such
    that if two elements $s$ and $t$ of $\mathcal{F}$ intersect
    inside $\H$, then the
     two associated generators $g_s$ and $ g_t$ must have the same fixed 
     point, which then, in particular,  has to be a parabolic fixed point of $G$ of rank  
     at least
     $2$.

    \noindent Let us now first recall from \cite{SS} the construction of the relevant 
    coding map $T$ associated with $G$,  which maps  the radial 
    limit set $L_{r}(G)$ into itself. This
    construction parallels the construction of the well-known  Bowen--Series map
    (cf. \cite{BowenSeries}, \cite{Sta1}, \cite{Sta2}).
    
   \noindent For $\xi, \eta \in  L_r(G)$, 
    let $\gamma_{\xi,\eta}: \R \to  \H$ denote to the directed geodesic 
    from
    $\eta$ to $\xi$ such that  $\gamma_{\xi,\eta}$ intersects 
    the closure $\overline{F}$ of $F$ in $\H$, and 
    normalised such that $\gamma_{\xi,\eta}(0)$ is
    the summit of
     $\gamma_{\xi,\eta}$. 
    The
      exit time  $e_{\xi,\eta}$ is
     defined  by
    \begin{eqnarray*}
    e_{\xi,\eta}  := \sup\{ s : \gamma_{\xi,\eta}(s) \in \overline{F} \} .
    \end{eqnarray*}
    Since  $\xi, \eta \in L_r(G)$, we clearly  have that
    $|e_{\xi,\eta}|<\infty$. 
    By  Poincar\'e's polyhedron theorem (cf. \cite{EP}),  we have that 
    the set  $\mathcal{F}$  carries an involution $\mathcal{F} \to
    \mathcal{F}$,  given by $s \mapsto s'$ and $ s''=s$. In 
    particular, for each $s \in
    \mathcal{F}$ there is a unique {face--pairing} transformation
    $g_s\in G$ such
    that  $g_s(\overline{F})\cap \overline{F} = s' $.
    We then let
	\[  \mathcal{L}_{r}(G) := \{(\xi,\eta) : \xi,\eta \in  L_r(G), \xi \neq 
	\eta \mbox{ and }\exists\,t\in \R\;:\; \gamma_{\xi,\eta}(t) 
	\in \overline{F}\},\]
    and define the map $S :\mathcal{L}_{r}(G)  
    \to \mathcal{L}_{r}(G) $,  for all  $( 
    \xi,\eta)\in \mathcal{L}_{r}(G) $
     such that $  \gamma_{\xi, \eta}(  e_{\xi,\eta}) \in s$, for some
    $ s \in  \mathcal{F}$, by
    \[ S(\xi,\eta) := (g_s (\xi),g_s (\eta)).  \]
    In order to
    show that the map $S$ admits a Markov partition, we
    introduce the
    following collection of subsets of the boundary $\partial \H$ of                                          %
$\H$. For $s \in \mathcal{F}$,
    let $A_s$ refer to the open hyperbolic halfspace
    for which $F \subset \H \setminus A_s$
    and $s \subset \partial A_s$. 
    Also, let $\Pi:\H \to \partial \H$ denote the 
  shadow-projection given by $\Pi(A):=\{\xi \in \partial \H: 
  \sigma_{\xi} \cap A \neq \emptyset \}$, where $\sigma_{\xi}$ denotes the ray from $0$ to $\xi$.    
    Then  the projections $a_s$ of the side
    $s$ to $ \partial \H$ is given by
    \[ a_s := \mbox{Int}(\Pi(A_s)).  \]
    If $G$ has exclusively parabolic fixed points 
    of rank $1$, then $a_s\cap
    a_t=\emptyset$, for all $ s,t \in \mathcal{F}, s\neq t$. Hence, by convexity
    of $F$,
    we have  $\gamma_{\xi,\eta}(e_{\xi,\eta}) \in s$ if and only if $\xi \in
    a_s$. In other words, $S(\xi, \eta) = (g_s \xi, g_s \eta)$ for all $\xi \in
    a_s$. This immediately gives that the projection map 
    $\pi:(\xi,\eta) \mapsto \xi$ onto the first coordinate
    of $\mathcal{L}_{r}(G) $ leads to a canonical factor $T$ of $S$,
    that is, we obtain the map
  \[  T: L_r(G) \to L_r(G), \mbox{ given
   by } T\arrowvert_{a_s \cap L_r(G)} := g_s.\]
 Clearly, $T$ satisfies $\pi \circ S = T\circ \pi$. Since
    $T(a_s) = g_s(a_s)= \mbox{Int}(\partial \H \setminus a_{s'})$, it follows that $T$  is 
    a non--invertible Markov
    map with respect to  the  partition $\{ a_s \cap L_r(G) : s \in
    \mathcal{F}  \}$. 

    If there are parabolic fixed points of rank greater than $1$, then, 
    a priori, $S$ does not have a canonical
     factor. In this  situation, the idea is to construct an 
     invertible  Markov map $\tilde{S}$ which is isomorphic
     to $S$ and which has a canonical factor. 
     This can be  achieved by   
 introducing a certain rule on the set  of faces associated with 
  the parabolic fixed points of rank greater than 1, which keeps 
 track 
    of the geodesic movement within these cusps. 
This then permits to define a
  coding map also in this higher 
  rank case, and, for ease of notation,  this map will also
  be denoted by $T$.
  (For further details we refer to \cite{SS}, where this construction 
  is given for $G$ having parabolic fixed points of rank $2$ and
  acting on $3$-dimensional hyperbolic space; the general case 
  follows from a straight forward adaptation of this construction.)
 For this so obtained coding map $T$ we then have the following result.
  \begin{prop}[{\cite[Proposition 2, Proposition 3]{SS}}]   
       The map $T$ is a topologically mixing Markov map with respect to the
    partition  generated by $\{ a_s \cap L_r(G): s \in
    \mathcal{F}\}$. 
    Moreover, the map $\tilde S$ is the natural
    extension of $T$.
  \end{prop} 
    \subsection{Patterson measure theory}     
    In order to introduce the  $T$-invariant measure 
    on  $ L(G)$ relevant for us here, 
    let us first briefly recall  some of the highlights in connection 
    with  the Patterson measure and the Patterson--Sullivan measure
    (for detailed discussions of these measures
     we refer to \cite{Patterson}, \cite{Su}, \cite{Su3},
     \cite{Nicholls}, \cite{SV}, see also \cite{DS} in these Proceedings). By now it is folklore that, given some sequence $(s_{n})$ of positive reals which tends to 
     $\delta$ from above, the Patterson measure 
     $m_{\delta}$ is a probability measure 
     supported on $L(G)$,   given by a weak accumulation point  of the sequence of measures
    \[ \left((\P_{\infty}(0,0,\delta_{n}))^{-1} \sum_{g\in G} 
    \e^{-\delta_{n}d(0,g(0))} \1_{g(0)}\right) .\]
  For geometrically finite Kleinian groups, and therefore, in 
  particular, for the type 
  of groups considered in this note, it is well known  that the so 
  obtained limit measure is non-atomic and does not depend on the particular choice of the
  sequence $(s_{n})$. Hence, in particular, $m_{\delta}$ is unique. 
 Moreover,  we have that $m_{\delta}$ is 
 $\delta$-conformal, that is, 
    for all $g\in G$ and $\xi \in
     L(G)$, we have 
  \[ \frac{d \, (m_{\delta}\circ g)}{d\, m_{\delta}} (\xi) = 
    \left(\frac{1-|g(0)|^{2}}{|\xi-g^{-1}(0)|^{2}}\right)^\delta.  \]
    This
    $\delta$-conformality is one of the key properties  of $m_{\delta}$,
    and for geome\-trically finite Kleinian groups
   it has the following, very useful geometrisation.  For this, 
     let
    $\xi_t$ denote the unique point on the ray $\sigma_{\xi}$  
    such that the
    hyperbolic distance between $0$
    and $\xi_{t}$  is equal to $t$, for arbitrary $\xi \in L(G)$ and $t>0$.  
    Also, let $B_{c}
    (\xi_t)\subset \H$ denote the
  $(N+1)$-dimensional hyperbolic disc  centred at $\xi_{t}$ of 
  hyperbolic radius $c>0$.
   Moreover, if $\xi_{t}$ 
    lies in one of the cusps  associated with the parabolic fixed 
    points of $G$,  we let
    $r(\xi_t)$ denote the rank of the parabolic fixed point associated 
    with that cusp, otherwise,  we put $r(\xi_t)$ equal
    to $\delta$.  We then have the following generalisation of Sullivan's 
    shadow lemma, where ``$\diam$'' denotes the Euclidean diameter 
    in $\partial \H$.
    \begin{prop}[\cite{SulivanNewOld},
    \cite{SV}]\label{measure}
    For  fixed, sufficiently large $c>0$, and for all $\xi \in L(G)$ 
    and $t >0$,  we have
     \[
    m_{\delta} ( \Pi(B_{c}( \xi_t ) )) \asymp 
    \left(\diam( \Pi(B_{c}( \xi_t ) ))\right)^{\delta} \cdot 
    \e^{ (
    r( \xi_t)-\delta)d(\xi_{t},G(0))}.  \]
    \end{prop} 
   \noindent A further strength of the Patterson measure in the geometrically 
   finite situation lies
   in the fact that it gives rise to
  a measure  $\tilde{m}_{\delta}$  on $(L(G) 
    \times L(G))\setminus \{\mbox{diag.}\}$, which is ergodic with 
    respect to the action of $G$ on $(L(G) 
    \times L(G))\setminus \{\mbox{diag.}\}$, given by 
    $g((\xi,\eta))=(g(\xi),g(\eta))$.  This measure is usually called 
    the  Patterson--Sullivan  measure, and it is given by 
    \[   d\tilde{m}_{\delta}(\xi,\eta) := \frac{d m_{\delta}(\xi) d 
    m_{\delta}(\eta)}{|\xi -\eta|^{2 \delta}} .      \]
    The (first) marginal measure of the Patterson--Sullivan measure 
    then defines the measure $\mu_{\delta}$ on $L(G)$, given by 
     \[   \mu_{\delta} := \tilde{m}_{\delta}\circ\pi^{-1}.\]
  The advantage of the measure  $\mu_{\delta}$ is that it is suitable 
  for non-trivial applications 
    of  certain results from infinite ergodic theory. 
    In fact,  for the system $\left(L(G),T,\mu_{\delta}\right)$  the following
       results have been 
       obtained in
       \cite{SS}.
    \begin{prop}\label{prop2}  The map  $T$ is conservative and ergodic 
	with respect 
    to the $T$-invariant, $\sigma$-finite measure $\mu_{\delta}$, and 
      $\mu_{\delta} $ is infinite if and only if $\delta \leq (r_{\maxx}
      +1)/2 $. Moreover,  if $G$ has  parabolic fixed points of rank less 
      than $r_{\maxx}$, then
      $\mu_{\delta}$
      gives  finite mass to   small  neighbourhoods around these fixed 
      points.
    \end{prop}
    \subsection{Infinite ergodic theory} In this section we 
    summarise some of the infinite ergodic theoretical properties
    of the system $\left(L(G),T,\mu_{\delta}\right)$. For further 
    details we refer to 
    \cite{SS}.
    
  \noindent  Recall that we always assume that $G$ is  a geometrically finite, 
    essentially free, zonal Kleinian group, and note 
    that for our purposes here we only have to consider the parabolic 
    subgroups of maximal rank, since, by Proposition \ref{prop2},
    $\mu_{\delta}$ gives infinite 
    measure to arbitrary small  neighbourhoods of a fixed point of a 
    parabolic generator  of $G$ 
    {\em only} if 
    the  parabolic fixed point is of maximal rank $r_{\maxx}$.
    Then define \[ \mathcal{D}_{0}:=  \bigcap_{ \gamma \mbox{ {\tiny a generator of} } \Gamma_{i} 
     \atop i=1,\ldots,v ;
     r_{i}=r_{\maxxx} }(\H \setminus \mbox{Cl}_{\H}(A_{\gamma \circ 
     \gamma})),\]
     and let
 \[ \mathcal{D}:=L_{r}(G) \cap \Pi(\mathcal{D}_{0}).\]
 Recall that the induced transformation $T_{\mathcal{D}}$ on 
 $\mathcal{D}$
 is defined by $T_{\mathcal{D}}(\xi):=T^{\rho (\xi)}(\xi)$,
 where $\rho$ denotes the return time function,
 given by $\rho(\xi):=\min \{n\in \N: T^{n}(\xi) \in
 \mathcal{D}\}$. 
One then considers the induced system 
 $\left(\mathcal{D},T_{\mathcal{D}}, \mu_{\delta,\mathcal{D}}\right) $,
 where  $\mu_{\delta,\mathcal{D}}$ denotes the restriction of
 $\mu_{\delta}$ to $\mathcal{D}$. Using standard techniques from ergodic 
 theory, for this induced system the following result was obtained in 
 \cite{SS}. Here,  $b_{n} \ll c_{n}$ means that  
$(b_{n}/c_{n})$
is uniformly bounded away from   infinity.
\begin{fact*}[\cite{SS}]
 The map $T_\mathcal{D}$ has the 
 Gibbs--Markov property with respect to the measure
    $\mu_{\delta,\mathcal{D}}$. That is, there exists $c \in (0,1)$ such that for
    arbitrary  cylinders 
 $[\omega_1]$ of length $n$ and $[\omega_2]$ of
 length $m$ such that $[\omega_2] \subset 
 T^n_\mathcal{D} ([\omega_1])$, 
 we have for  $\mu_{\delta,\mathcal{D}}$-almost
 every pair $\eta,\xi \in [\omega_2]$,   
\[\left|  \log \frac{d \mu_{\delta,\mathcal{D}}\circ 
 T^{-n}_{\mathcal{D},
 \omega_1}}{d\mu_{\delta,\mathcal{D}}}(\xi) - \log\frac{d 
 \mu_{\delta,\mathcal{D}}\circ
     T^{-n}_{\mathcal{D},\omega_1}}{d\mu_{\delta,\mathcal{D}}}(\eta) 
     \right| 
     \ll c^m ,\]
 where $ T^{-n}_{\mathcal{D},\omega_1}$ denotes the inverse 
 branch of $T_\mathcal{D}^n$
 mapping $T^{n}_{\mathcal{D}}([\omega_1]) $ bijectively to $[\omega_1]$.
\end{fact*}
\noindent Let   $\widehat T_\mathcal{D}$ denotes the dual operator of $T_\mathcal{D}$, given by
\[
 \mu_{\delta, \mathcal{D}}(f\cdot g\circ T)= \mu_{\delta, 
 \mathcal{D}}(\widehat{T}_\mathcal{D}(f)\cdot g), 
 \mbox{ for all }f\in L^1( \mu_{\delta, 
 \mathcal{D}}),g\in L^\infty ( \mu_{\delta, \mathcal{D}}).
\] 
The Gibbs--Markov property of $T_\mathcal{D}$  then allows to employ the following
 chain of implications (cf. \cite{Aa2},
\cite{ADU}):

{\it $T_\mathcal{D}$ has the 
Gibbs--Markov property with respect to 
$\mu_{\delta,\mathcal{D}}$.
%\vspace{-.1cm}
\begin{itemize}
\item[$\implies$] 
There
  exists $c_{0}\in (0,1)$ such that, for
  all $f \in L^1(\mu_{\delta,\mathcal{D}})$ and $n \in \N$, we have
  \[ \left\| {\widehat T_\mathcal{D}}^n f - \mu_{\delta, \mathcal{D}}( 
  f) \right\|_L \ll c_{0}^n \|f\|_L
  . \]
  (Here, $\| \cdot \|_L$ refers to the Lipschitz norm
  (cf. \cite{ADU}, p. 541).)  
\item[$\implies$] $T_\mathcal{D}$ is continued fraction mixing 
(cf. \cite{ADU}, p. 500).
\item[$\implies$] The set $\mathcal{D}$ is a Darling--Kac set 
for $T$. That is,   there  exists a sequence 
$(\nu_n)$ (called the {\em return
  sequence} of $T$) such that 
  \[ \frac{1}{\nu_n} \sum_{i=0}^{n-1} \widehat{ T}^i \1_\mathcal{D}(\xi) 
  \to \mu_{\delta}(\mathcal{D}), \quad
  \mbox{uniformly for } \mu_{\delta} \mbox{-almost every } \xi 
  \in \mathcal{D}.\]
\end{itemize} }
\noindent Finally, let us also remark that  the growth rate of the sequence $(\nu_n)$
can be determined explicitly as follows.  
	  Recall from \cite[Section 3.8]{Aa2}
	  that  
	  the wandering rate of the Darling--Kac set  $\mathcal{D}$ is defined
	  by the sequence 
	  $\left(w_{n}\right)$,
	  which is given, for each $n \in \N$, by  \[
	  w_{n}:=
	  \mu_{\delta}\left(\bigcup_{k=1}^{n}
	  T^{-(k-1)}(\mathcal{D})\right).\]
%Since $(L_r(G), T, \mu_{\delta})$ has a Darling-Kac set, we can apply 
An application of  \cite[Proposition 3.8.7]{Aa2}  gives that the return sequence and  the wandering rate are related through 
	   \[\nu_{n}\cdot 
	   w_{n}\sim \frac{n}{\Gamma(1+\beta) \Gamma(2-\beta)},\]
where $\beta:=\max \{0,1+r_{\maxx} -2 \delta\}$ coincides with the index of 
variation of the regularly varying sequence $(w_{n})$.	   
Hence, we are left with to determine the wandering rate. But this
has been done in \cite[Theorem 1]{SS}, where it was shown 
that
\[ w_n \asymp  \left\{ \begin{array}[h]{l @{$\quad\mbox{  for  }\quad$} l}
 n^{ r_{\maxxx} -2 \delta +1} & \delta <( r_{\maxxx} +1)/2 \\
 \log n & \delta =( r_{\maxxx} +1)/2 \\
	1 & \delta >( r_{\maxxx} +1)/2 .
\end{array} \right.  \]
Hence, by combining these observations, it follows that
\[ \nu_n \asymp \left\{
  \begin{array}[h]{l @{\, \,\,  \hbox{  for  } \, \, \, }l}
n^{2\delta - r_{\maxxx}} & \delta < ( r_{\maxxx}
 +1)/2\\
n/\log n &  \delta = (r_{\maxxx}  +1 )/2\\
n &  \delta > ( r_{\maxxx}  +1)/2 .
  \end{array}
\right.
\]

    \section{Proof of the Theorem \ref{main}}
   \noindent As  we have seen in the previous section, we have  that 
   the set $\mathcal{D}:=L_{r}(G) \cap \Pi(\mathcal{D}_{0})$
     is a 
    Darling--Kac set.   Combining this with Proposition \ref{measure}, 
    Proposition \ref{prop2}, and the fact that the Patterson measure $m_{\delta}$
    and its $T$--invariant version $\mu_{\delta}$ are comparable on $\mathcal{D}$, one obtains
    \begin{eqnarray*}
	\frac{1}{\nu_{n}}\sum_{g(w)\in \mathcal{D}_{0} \atop
	|g|\leq n}\e^{-\delta d(z,g(w))} 
	& \asymp & \frac{1}{\nu_{n}}\sum_{k=0}^{n}m_{\delta}
	(\mathcal{D}\cap 
	T^{-k}(\mathcal{D}))\asymp\frac{1}{\nu_{n}}
	\sum_{k=0}^{n}
	\mu_{\delta}(\mathcal{D}\cap
	T^{-k}(\mathcal{D}))\\
	 & = & \frac{1}{\nu_{n}}\sum_{k=0}^{n}
	 \mu_{\delta}(\1_{\mathcal{D}}
	 \cdot\hat{T}^{k}\1_{\mathcal{D}})=\mu_{\delta}
	 (\1_{\mathcal{D}}
	 \cdot\frac{1}{\nu_{n}}\sum_{k=0}^{n}\hat{T}^{k}
	 \1_{\mathcal{D}})\\
	 & \sim & \left(\mu_{\delta}(\1_{\mathcal{D}})
	 \right)^{2}.
	     \end{eqnarray*}
 Since $\mu_{\delta}(\1_{\mathcal{D}})\asymp 1$, it follows that
    \[  \sum_{g(w)\in \mathcal{D}_{0}\atop
    |g|\leq n} \e^{-\delta d(z,g(w))} \asymp \nu_{n} 
    .\]
    To extend this estimate to the full $G$-orbit of $w$, let
    \[Q_{i}:=  \bigcap_{ \gamma \mbox{ {\tiny a generator of} } 
    \Gamma_{i} }(\H \setminus \mbox{Cl}_{\H}(A_{\gamma}))\] 
    denote the fundamental domain for the action of $\Gamma_i$ on 
    $\H$, for each $i \in \{1,\ldots, v\}$. Clearly,  we can assume, without loss of generality,
    that $z$ and $w$ are 
    contained in each of the fundamental domains $Q_{i}$. For every $\gamma \in \Gamma_i$
    such that $|\gamma|=k$, for some $1<k\leq n$, we then have, 
    with the convention $\nu_0:=1$,
       \[
    \sum_{g(w) \in \gamma(Q_i) \atop |g|\leq n}\e^{-\delta d(z,g(w))}\asymp k^{-2\delta}\nu_{n-k}.\]
 Also, note that 
 \[ \card \{\gamma \in \Gamma_i:|\gamma| = k \} \asymp k^{r_i -1}.\]
 Combining these observations,
it follows that
\begin{eqnarray*}
    \P_{n}(\delta,z,w) & \asymp &  \sum_{g(w)\in \mathcal{D}_{0}\atop
    |g|\leq n} \e^{-\delta d(z,g(w))} +  \sum_{i=1,..,v \atop
    r_i = r_{\maxxx}} \sum_{\gamma \in 
    \Gamma_i \atop
    |\gamma|\geq 2} \sum_{g(w) \in \gamma(Q_i) \atop |g|\leq n}
    \e^{-\delta d(z,g(w))}\\
    &
     \asymp & \nu_n	+  \sum_{i=1,..,v \atop
    r_i = r_{\maxxx}} \sum_{k=2}^{n} \sum_{\gamma \in 
    \Gamma_i \atop |\gamma| = k}  \sum_{g(w) \in 
    \gamma(Q_i) \atop |g| \leq n} \e^{-\delta d(z,g(w))}\\
   \\
    & \asymp & \nu_n	+  \sum_{k=2}^{n}  k^{r_{\maxxx} -1}  k^{-2\delta}\nu_{n-k}.
\end{eqnarray*}
In order to finish the proof, recall that, 
by a result of Beardon \cite{Beardon}, one has that  $\delta>r_{\maxx}/2$.
Therefore,  there exists $\kappa=\kappa(G)>0$  such that
$\delta>r_{\maxx}/2+\kappa$. Moreover, note that  we can assume,
without loss of generality, that
    $\left(\nu_{n}\right)$ is increasing. Using these observations, it 
    now follows that, on the one hand, \[
   \sum_{k=2}^{n}\nu_{n-k}k^{-2\delta+r_{\maxxx}-1}
   \ll \nu_{n}\sum_{k=2}^{n}k^{-2\delta+r_{\maxxx}-1}\ll
  \nu_{n}\sum_{k=1}^{n}k^{-1-2\kappa}\ll \nu_{n}.\]
  On the other hand, we clearly have that
  $\sum_{k=2}^{n}\nu_{n-k}k^{-2\delta+r_{\maxxx}-1} \gg \nu_{n-2}$.
  Combining these observations with the estimate for the asymptotic growth rate of the 
  return sequence $(\nu_{n})$, given in the previous section, the proof of 
  Theorem \ref{main} follows.
 
\singlespace

\end{document}